# A LONG RANGE DEPENDENCE STABLE PROCESS AND AN INFINITE VARIANCE BRANCHING SYSTEM[1]

By Tomasz Bojdecki, Luis G. Gorostiza and Anna Talarczyk

*University of Warsaw, Centro de Investigación y de Estudios Avanzados and University of Warsaw*

We prove a functional limit theorem for the rescaled occupation time fluctuations of a $(d,\alpha,\beta)$-branching particle system [particles moving in $\mathbb{R}^d$ according to a symmetric $\alpha$-stable Lévy process, branching law in the domain of attraction of a $(1+\beta)$-stable law, $0 < \beta < 1$, uniform Poisson initial state] in the case of intermediate dimensions, $\alpha/\beta < d < \alpha(1+\beta)/\beta$. The limit is a process of the form $K\lambda\xi$, where $K$ is a constant, $\lambda$ is the Lebesgue measure on $\mathbb{R}^d$, and $\xi = (\xi_t)_{t \geq 0}$ is a $(1+\beta)$-stable process which has long range dependence. For $\alpha < 2$, there are two long range dependence regimes, one for $\beta > d/(d+\alpha)$, which coincides with the case of finite variance branching ($\beta = 1$), and another one for $\beta \leq d/(d+\alpha)$, where the long range dependence depends on the value of $\beta$. The long range dependence is characterized by a dependence exponent $\kappa$ which describes the asymptotic behavior of the codifference of increments of $\xi$ on intervals far apart, and which is $d/\alpha$ for the first case (and for $\alpha = 2$) and $(1+\beta-d/(d+\alpha))d/\alpha$ for the second one. The convergence proofs use techniques of $\mathcal{S}'(\mathbb{R}^d)$-valued processes.

**1. Introduction.** We consider a so-called $(d,\alpha,\beta)$-branching particle system in $\mathbb{R}^d$ described as follows. Particles start off at time $t = 0$ from a Poisson random field with intensity measure $\lambda \equiv \lambda_d$ (Lebesgue measure), and they evolve independently, moving according to a standard symmetric $\alpha$-stable Lévy process ($0 < \alpha \leq 2$) and splitting at rate $V$ according to the branching

---

Received August 2005; revised May 2006.

[1]Supported in part by CONACyT Grant 45684-F (Mexico) and MEiN Grant 1P03A01129 (Poland).

law

$$p_k = \begin{cases} \dfrac{1}{1+\beta} \binom{1+\beta}{k} (-1)^k, & k = 0, 2, 3, \ldots, \\ 0, & k = 1, \end{cases}$$

$0 < \beta \leq 1$. This law is critical, for $\beta = 1$, it is binary branching, and for $\beta < 1$, it is a simple typical element of the domain of attraction of a stable law with exponent $1 + \beta$. Its generating function is

(1.1) $$s + \frac{1}{1+\beta}(1-s)^{1+\beta}, \qquad 0 < s < 1.$$

This branching particle system and its associated superprocess have been widely studied; some of the early results appear in [9, 10, 11, 14, 17, 19, 20] and references therein. In this paper we prove a functional limit theorem for the rescaled occupation time fluctuations of the system with $\beta < 1$ in the case of "intermediate dimensions," $\alpha/\beta < d < \alpha(1+\beta)/\beta$, where the limit process has long range dependence (Theorem 2.2), and we show that for $\alpha < 2$ there are two different long range dependence regimes, depending on whether $\beta$ is above or below the value $d/(d+\alpha)$ (Theorem 2.7).

Let us recall first the result for $\beta = 1$, where the branching law has finite variance. A functional limit theorem for the occupation time fluctuations of the system was proven in [6] when the limit process has long range dependence, which occurs for the intermediate dimensions

(1.2) $$\alpha < d < 2\alpha.$$

The limit process is of the form $C\lambda\zeta$, where $C$ is a constant and $\zeta = (\zeta_t)_{t \geq 0}$ is a real valued, continuous, self-similar, long range dependence Gaussian process, called *sub-fractional Brownian motion*, whose covariance function is

$$s^h + t^h - \tfrac{1}{2}[(s+t)^h + |s-t|^h],$$

where $h = 3 - d/\alpha$. The long range dependence of $\zeta$ is characterized by the behavior of the covariance of increments on intervals separated distance $T$, which decays like $T^{-d/\alpha}$ as $T \to \infty$. Properties of $\zeta$ are studied in [5]. The reason for the name sub-fractional Brownian motion is the fact that the increments of $\zeta$ on nonoverlapping intervals are more weakly correlated than those of fractional Brownian motion, whose covariance function is $\tfrac{1}{2}[s^h + t^h - |s-t|^h]$, and their covariance decays faster as the distance between the intervals tends to $\infty$; in this sense sub-fractional Brownian motion is intermediate between Brownian motion and fractional Brownian motion.

The Gaussian property of the limit process in [6] is due to the finiteness of the variance of the binary branching law. More general critical, finite variance branching laws lead to essentially the same limit process [23]. A substantially different and more interesting situation occurs with infinite variance branching, $\beta < 1$. There are new technical problems in the proofs and the results reveal new properties of the $(d, \alpha, \beta)$-branching particle system. Condition (1.2) for long range dependence is now replaced by

$$\text{(1.3)} \qquad \frac{\alpha}{\beta} < d < \frac{\alpha(1+\beta)}{\beta}$$

(note that $\alpha/\beta < d$ is the condition for the system to converge for large time toward an equilibrium state which has intensity $\lambda$; for $d \leq \alpha/\beta$, the system goes to local extinction [19]). The occupation time fluctuation limit process for $\beta < 1$ resembles the one for $\beta = 1$ in that it has a simple spatial structure and a complicated temporal one. It has the form $K\lambda\xi$, where $K$ is a constant and $\xi = (\xi_t)_{t \geq 0}$ is a continuous, self-similar, long range dependence $(1+\beta)$-stable process, which may be called a *sub-fractional stable process* by analogy with the case $\beta = 1$ [see Remark 2.4(a)]. The convergence takes place in the space of continuous functions $C([0, \tau], \mathcal{S}'(\mathbb{R}^d))$ for any $\tau > 0$, where $\mathcal{S}'(\mathbb{R}^d)$ is the space of tempered distributions, dual of the space $\mathcal{S}(\mathbb{R}^d)$ of smooth rapidly decreasing functions. We stress that convergence methods for $\mathcal{S}'(\mathbb{R}^d)$-valued processes play a fundamental role in this paper due to the convenient topology of $\mathcal{S}'(\mathbb{R}^d)$ (although here all infinite dimensional processes are measure-valued).

Concerning functional convergence of the occupation time fluctuation process, in the finite variance case the tightness proof in [6] employed standard methods based on moment estimates. For $\beta < 1$, there are no moments of orders $\geq 1 + \beta$ and a more delicate approach is needed. The general scheme for the identification of a unique limit is similar to that in [6], which involves a space-time random field method introduced in [4]. This approach is simpler in the present case than proving convergence of finite dimensional distributions, but additional technical work is needed to handle $\beta < 1$, due to the fact that the Fourier transform method that was widely used in [6] is not applicable in most cases in this paper.

The long range dependence of the process $\xi$ is characterized by means of the asymptotic behavior of the codifference of increments on intervals distance $T$ apart as $T \to \infty$, and this is given in terms of a *dependence exponent* $\kappa$ (Definition 2.5). Regarding codifference, see [25]. For $\alpha = 2$, the dependence exponent is $\kappa = d/\alpha$, and the codifference decays at the same rate as the covariance in the finite variance case ($\beta = 1$), that is, $T^{-d/\alpha}$ [5]. For $\alpha < 2$, there are two long range dependence regimes which are separated by the value $\beta = d/(d + \alpha)$:

(1) For $d/(d+\alpha) < \beta < 1$, the dependence exponent is $\kappa = d/\alpha$.
(2) For $\beta \leq d/(d+\alpha)$, the dependence exponent depends on $\beta$ as follows: $\kappa = (1 + \beta - d/(d+\alpha))d/\alpha$.

Our methods are mainly analytic and they have technical advantages, but "physical" explanations of the results in terms of the particle system are not apparent from them. The simple spatial structure of the limit process (Lebesgue measure) arises already in occupation time fluctuations of a Poisson system of independent Brownian motions in dimension $d = 1$ [15] [Theorem 0.4(i)], and the same thing happens with critical binary branching two dimensions higher, $d = 3$ [9] (Theorem 1), which is the special case of (1.3) for $\alpha = 2$, $\beta = 1$ (see also [12] for a more general setup which also shows the shift in dimensions from nonbranching to branching; [9] and [12] have only single time limits). The cause of the long range dependence with $\beta = 1$ was attributed intuitively in [6] to "clan recurrence." A clan is a family of particles with a common ancestor. Clan recurrence means that any ball is visited by clans infinitely often and at arbitrarily large times, each visit adding a random amount to the occupation time of the ball. This phenomenon was investigated in [29] for the branching particle system in equilibrium with $\alpha = 2$, $\beta = 1$ (and a partial result for $\beta < 1$); it was shown that clan recurrence holds in dimensions $d = 3, 4$. The $(d, \alpha, \beta)$-branching system with initial Poisson condition is not in equilibrium. However, based on [6, 7, 8, 29], and the results in the present paper, the natural conjectures are that clan recurrence holds for the general $(d, \alpha, \beta)$-branching system if and only if $\alpha/\beta < d \leq \alpha(1 + \beta)/\beta$, with equilibrium or some other initial conditions (including Poisson), and that this is the cause of the long range dependence of the occupation time fluctuation limit process. But at the borderline $d = \alpha(1 + \beta)/\beta$, the long range dependence disappears [7, 8]. We would like to have some intuitive explanation for the two long range dependence regimes separated by $\beta = d/(d + \alpha)$ (whereas no such phenomenon occurs in the Brownian case, $\alpha = 2$), but this still remains an enigma to us. A related curious fact is that the dependence exponent $\kappa$ converges to $d/(d + \alpha)$ as $\beta \searrow 0$.

We mention some other results on occupation times. As already stated, for $d \leq \alpha/\beta$, the system becomes extinct locally. However, for $d = \alpha/\beta$, there is a functional ergodic theorem [30]. With $\beta = 1$, there are functional limit theorems for the fluctuations in dimensions $d \geq 2\alpha$, where there is no long range dependence [7]. In [3] the fluctuations of the occupation time of the origin are studied for a critical binary branching random walk on the $d$-dimensional lattice, $d \geq 3$; the convergence results are parallel to those in [6] and [7], but the proofs are quite different. [3] treats also the case of the branching random walk in equilibrium (see also [5] for the case $\beta = 1$ in equilibrium, where only covariance calculations were done). The results of [6] and [7] are extended in [23] for general critical, finite variance branching, with both equilibrium and Poisson initial conditions. [6] and [7] contain references to other relevant papers, among them those that awakened our interest in this subject, [9], [15] and [20] (although they do not refer to long range dependence). In [8] we study the occupation time fluctuations of the $(d, \alpha, \beta)$-branching particle system with $\beta < 1$ and Poisson initial condition in dimensions $d \geq \alpha(1 + \beta)/\beta$, where there is no long range dependence. For the "critical" dimension, $d = \alpha(1 + \beta)/\beta$, the limit process has the form $K\lambda\eta$, where $\eta = (\eta_t)_{t \geq 0}$ is a real $(1 + \beta)$-stable process with stationary independent increments. For "large" dimensions, $d > \alpha(1 + \beta)/\beta$, the limit is an $\mathcal{S}'(\mathbb{R}^d)$-valued (not measure-valued) $(1 + \beta)$-stable process with stationary independent increments. The size of the occupation time fluctuations is given by the norming for convergence. If $T$ is the time scaling parameter and $F_T$ denotes the norming [see (1.4) below], $F_T$ is $T^{(2+\beta-(d/\alpha)\beta)/(1+\beta)}$ for $\alpha/\beta < d < \alpha(1+\beta)/\beta$, $(T \log T)^{1/(1+\beta)}$ for $d = \alpha(1+\beta)/\beta$, and $T^{1/(1+\beta)}$ for $d > \alpha(1+\beta)/\beta$. This illustrates the phenomena that atypical normings are associated with long range dependence, and higher order fluctuations (involving a logarithmic factor) arise at the borderline between different behaviors. The norming for high dimensions is the one for the classical central limit theorem. Several unsolved questions of "physical" interpretation are also brought up in [8], including the ones mentioned above.

Long range dependence is now an area of intensive research due to its mathematical appeal and its manyfold applications (see, e.g., [16]). In particular, there are other long range dependence, infinite variance processes, for example, [21]. Our interest in the subject was inspired by the appearance of long range dependence in occupation time fluctuations of particle systems, specially branching systems and related superprocesses. Other types of long range dependence Gaussian processes connected with branching systems with immigration are presented in [18] (without functional convergence proofs).

We now give some definitions and notation.

For the $(d, \alpha, \beta)$-branching particle system (with $\beta < 1$), let $(N_t)_{t \geq 0}$ denote the empirical measure process, that is, $N_t(A)$ is the number of particles in the set $A \subset \mathbb{R}^d$ at time $t$. Thus, $N_0$ is a Poisson random measure with intensity $\lambda$. The rescaled occupation time fluctuation process is defined by

$$(1.4) \quad X_T(t) = \frac{1}{F_T} \int_0^{Tt} (N_s - \lambda) \, ds = \frac{T}{F_T} \int_0^t (N_{Ts} - \lambda) \, ds, \qquad t \geq 0,$$

where $F_T$ is a norming to be determined, and $T$ is the scaling parameter which accelerates the time and will tend to $\infty$. Note that $EN_s = \lambda$ for all $s$, due to the initial Poisson condition, the criticality of the branching and the $\alpha$-stable motion.

Constants are written $C, C_1$, and so forth, with possible dependencies in parenthesis. $\langle \cdot, \cdot \rangle$ denotes pairing of spaces in duality [e.g., $\mathcal{S}'(\mathbb{R}^k)$ and $\mathcal{S}(\mathbb{R}^k)$]. $\Rightarrow$ stands for weak convergence.

Section 2 contains the results, and Sections 3 and 4 the proofs.

**2. Results.** We start by introducing the process that plays a fundamental role in the paper. Let $M$ be the independently scattered $(1+\beta)$-stable measure on $\mathbb{R}^{d+1}$ with control measure $\lambda_{d+1}$ (Lebesgue measure) and skewness intensity 1, that is, for each $A \in \mathcal{B}(\mathbb{R}^{d+1})$ such that $0 < \lambda_{d+1}(A) < \infty$, $M(A)$ is a $(1+\beta)$-stable random variable with characteristic function

$$\exp\left\{-\lambda_{d+1}(A)|z|^{1+\beta}\left(1 - i(\operatorname{sgn} z)\tan\frac{\pi}{2}(1+\beta)\right)\right\}, \qquad z \in \mathbb{R},$$

the values of $M$ are independent on disjoint sets, and $M$ is $\sigma$-additive a.s. (see [26], Definition 3.3.1). Let $p_t(x)$ denote the transition density of the symmetric $\alpha$-stable Lévy process in $\mathbb{R}^d$, and recall our assumption

$$(2.1) \qquad \frac{\alpha}{\beta} < d < \frac{\alpha(1+\beta)}{\beta}.$$

The process $\xi = (\xi_t)_{t \geq 0}$ is defined as follows.

DEFINITION 2.1. Let

$$(2.2) \qquad \xi_t = \int_{\mathbb{R}^{d+1}} \left(\mathbb{1}_{[0,t]}(r) \int_r^t p_{u-r}(x)\,du\right) M(dr\,dx), \qquad t \geq 0,$$

where the integral with respect to $M$ is understood in the sense of [26], (3.2)–(3.4).

By [26], existence of this process follows from the fact that

$$(2.3) \qquad \int_{\mathbb{R}^d} \int_0^t \left(\int_r^t p_{u-r}(x)\,du\right)^{1+\beta} dr\,dx < \infty.$$

It can be verified that under (2.1) this integral is indeed finite (see, e.g., [17], Lemma A.1).

The first main result is the following functional limit theorem for the process $X_T$ defined by (1.4).

THEOREM 2.2. *Assume* (2.1). *Let*

$$(2.4) \qquad F_T = T^{(2+\beta-(d/\alpha)\beta)/(1+\beta)}.$$

*Then*

$$X_T \Rightarrow K\lambda\xi$$

*in* $C([0,\tau], \mathcal{S}'(\mathbb{R}^d))$ *as* $T \to \infty$ *for any* $\tau > 0$, *where* $\xi$ *is the process defined by* (2.2) *and*

$$K = \left(-\frac{V}{1+\beta}\cos\frac{\pi}{2}(1+\beta)\right)^{1/(1+\beta)}.$$

In the next proposition we collect some basic properties of the process $\xi$.

PROPOSITION 2.3. (a) $\xi$ is $(1+\beta)$-stable, totally skewed to the right, with finite-dimensional distributions given by

$$
\begin{aligned}
&E\exp\{i(z_1\xi_{t_1}+\cdots+z_k\xi_{t_k})\}\\
&=\exp\Bigg\{-\int_{\mathbb{R}^{d+1}}\left|\sum_{j=1}^k z_j\mathbb{1}_{[0,t_j]}(r)\int_r^{t_j}p_{u-r}(x)\,du\right|^{1+\beta}\\
&\qquad\times\left[1-i\operatorname{sgn}\left(\sum_{j=1}^k z_j\mathbb{1}_{[0,t]}(r)\int_r^{t_j}p_{u-r}(x)\,du\right)\right.\\
&\qquad\qquad\left.\times\tan\frac{\pi}{2}(1+\beta)\right]dr\,dx\Bigg\},
\end{aligned}
$$
(2.5)

$0\leq t_1<\cdots<t_k, z_1,\ldots,z_k\in\mathbb{R}$.

(b) $\xi$ is self-similar with index $H=(2+\beta-\frac{d}{\alpha}\beta)/(1+\beta)$, that is,

$$(\xi_{at_1},\ldots,\xi_{at_k})\stackrel{d}{=}a^H(\xi_{t_1},\ldots,\xi_{t_k}),\qquad a>0.$$

(c) $\xi$ has continuous paths (more precisely, has a continuous version).

Property (a) follows immediately from the definition and [26] (Proposition 3.4.2). Property (b) can be easily derived from (2.5) using the self-similarity of $p_t$. Property (c) is a consequence of Theorem 2.2.

REMARK 2.4. (a) It is not hard to verify that if we put $\beta=1$ in (2.5), we obtain the finite-dimensional distributions of the sub-fractional Brownian motion (multiplied by a constant) considered in [5] and [6]. This, together with the fact that Theorem 2.2 is an analogue of Theorem 2.2 in [6], suggests giving the name *sub-fractional stable process* to $\xi$. See [5] concerning relationships between sub-fractional Brownian motion and fractional Brownian motion, which make the name more appropriate in that case.

(b) We think that a functional limit theorem also holds for the occupation time fluctuation process of the $(d,\alpha,\beta)$-system with $\beta<1$ and initial equilibrium state, but we have not endeavored to prove it. In this case we conjecture that the limit process is of the form $K\lambda\eta$, where $\eta=(\eta_t)_{t\in\mathbb{R}}$ is a self-similar, continuous, stable process with stationary increments, which should be a kind of fractional stable process. Moreover, by analogy with the covariance results for the case $\beta=1$ [5], we conjecture that the process $\xi$ defined by (2.2) has the same distribution as the process $(\eta_t+\eta_{-t})_{t\geq 0}$ (multiplied by a constant).

(c) The continuity of $\xi$ can also be derived from the results of [22] with some technical work.

The process $\xi$ does not have independent increments and the increments are not stationary. In order to investigate its long range dependence, we introduce the following general notion.

DEFINITION 2.5. Let $\eta$ be a real infinitely divisible process. For $0 \leq u < v < s < t, T > 0, z_1, z_2 \in \mathbb{R}$, let

$$
\begin{aligned}
D_T&(z_1, z_2; u, v, s, t) \\
&= |\log E e^{i(z_1(\eta_v - \eta_u) + z_2(\eta_{T+t} - \eta_{T+s}))} \\
&\quad - \log E e^{iz_1(\eta_v - \eta_u)} - \log E e^{iz_2(\eta_{T+t} - \eta_{T+s})}|.
\end{aligned}
\tag{2.6}
$$

We define the *dependence exponent* $\kappa$ of the process $\eta$ by

$$
\kappa = \inf_{z_1, z_2 \in \mathbb{R}} \inf_{0 \leq u < v < s < t} \sup\{\gamma > 0 : D_T(z_1, z_2; u, v, s, t) = o(T^{-\gamma}) \text{ as } T \to \infty\}.
\tag{2.7}
$$

REMARK 2.6.  (a) If $\eta$ has independent increments, then $\kappa = +\infty$.

(b) If $\eta$ is Gaussian, then

$$D_T(z_1, z_2; u, v, s, t) = |z_1 z_2 \operatorname{Cov}(\eta_v - \eta_u, \eta_{T+t} - \eta_{T+s})|.$$

(c) $D_T$ is the modulus of the codifference of the random variables $z_1(\eta_v - \eta_u)$ and $-z_2(\eta_{T+t} - \eta_{T+s})$, as defined in [25] (see also [26] for symmetric stable case).

The second main result is the following theorem on the long range dependence of the process $\xi$.

THEOREM 2.7. *The dependence exponent of the process $\xi$ defined by* (2.2) *is given by*

$$
\kappa = \begin{cases} \dfrac{d}{\alpha}, & \text{if either } \alpha = 2, \text{ or } \alpha < 2 \text{ and } \beta > \dfrac{d}{d+\alpha}, \\ \dfrac{d}{\alpha}\left(1 + \beta - \dfrac{d}{\alpha+d}\right), & \text{if } \alpha < 2 \text{ and } \beta \leq \dfrac{d}{d+\alpha}. \end{cases}
\tag{2.8}
$$

REMARK 2.8.  (a) Note that for $\beta > 1/\sqrt{2}$, we have $\beta > d/(d+\alpha)$, and for $\beta < (\sqrt{5} - 1)/2$, we have $\beta < d/(d+\alpha)$.

(b) As we shall see in the proof, in the case $\beta > d/(d+\alpha)$ the value of $\kappa$ gives the exact asymptotics of $D_T$, that is,

(2.9) $$D_T(z_1, z_2; u, v, s, t) = O(T^{-d/\alpha}) \quad \text{as } T \to \infty,$$

provided that $z_1 z_2 > 0$.

(c) The real valued limit process in Theorem 2.2 in [6] has the form $C(\alpha, d)\zeta$, where the process $\zeta$ (sub-fractional Brownian motion) depended only on $d/\alpha$. In the present case all relevant parameters related to the process $\xi$ depend only on $\beta$ and $d/\alpha$ [see (2.1), (2.4), Proposition 2.3(b) and (2.8)]. A natural question is whether $\xi$ also has a form $C(\alpha, d, \beta)\xi'$, where the distribution of the process $\xi'$ depends only on $\beta$ and $d/\alpha$. We have not been able to answer this question.

(d) The standard symmetric $\alpha$-stable Lévy process on $\mathbb{R}^d$ is transient for $d > \alpha$, and its *degree of transience*, defined as

$$\gamma = \sup\{\theta > 0 : EL^\theta < \infty\},$$

where $L$ is the last exit time from an open unit ball centered at the origin, is given by

$$\gamma = \frac{d}{\alpha} - 1$$

(see [13] and [27]). In the $(\gamma, \beta)$-plane the regions corresponding to the two long range dependence regimes of the process $\xi$ are given as follows:

$$\kappa = \gamma + 1 \quad \text{for } 0 < \gamma < \sqrt{2} \text{ and } \max\left\{\frac{1}{\gamma+1}, \frac{\gamma+1}{\gamma+2}\right\} < \beta < \min\left\{1, \frac{1}{\gamma}\right\},$$

$$\kappa = (\gamma+1)\left(1 + \beta - \frac{\gamma+1}{\gamma+2}\right)$$

$$\text{for } \gamma > \frac{\sqrt{5}-1}{2} \text{ and } \frac{1}{\gamma+1} < \beta < \min\left\{\frac{1}{\gamma}, \frac{\gamma+1}{\gamma+2}\right\},$$

and $\kappa = \gamma + 1$ on the separating curve $\beta = \frac{\gamma+1}{\gamma+2}$, $\frac{\sqrt{5}-1}{2} < \gamma < \sqrt{2}$.

**3. Proof of Theorem 2.2.** Without loss of generality we assume $\tau = 1$. To start, we gather some technical facts which will be used in the proof several times.

Recall that $p_t(\cdot)$ has characteristic function $e^{-t|z|^\alpha}$. We denote by $\mathcal{T}_t$ the corresponding semigroup, that is, $\mathcal{T}_t f = p_t * f$. It is well known that (self-similarity)

(3.1) $$p_t(x) = t^{-d/\alpha} p_1(x t^{-1/\alpha}),$$

and for $\alpha < 2$,

(3.2) $$\frac{c_1}{1 + |x|^{d+\alpha}} \leq p_1(x) \leq \frac{c_2}{1 + |x|^{d+\alpha}}$$

for some positive constants $c_1$ and $c_2$. The upper bound obviously holds also for $\alpha = 2$.

We will use the following slightly stronger version of (2.3), obtained by the same argument:

$$\text{(3.3)} \qquad \int_{\mathbb{R}^d} \left( \int_0^1 p_u(x)\, du \right)^{1+\beta} dx < \infty,$$

provided that $d < \alpha(1+\beta)/\beta$ [see (2.1)].

We will also use the following two elementary estimates:

$$\text{(3.4)} \qquad \begin{aligned} 0 &\leq (a+b)^{1+\beta} - a^{1+\beta} - b^{1+\beta} \\ &\leq (1+\beta) a^\delta b^{1+\beta-\delta}, \qquad a, b \geq 0, \beta \leq \delta \leq 1, \end{aligned}$$

$$\text{(3.5)} \qquad (a+b)^{1+\beta} - a^{1+\beta} - b^{1+\beta} \geq \beta b^\beta a, \qquad b \geq a \geq 0.$$

We will employ the space-time random field approach [4], which consists of investigating weak convergence of the $\mathcal{S}'(\mathbb{R}^{d+1})$-random variable $\widetilde{X}_T$ associated with the process $X_T$, defined by

$$\text{(3.6)} \qquad \langle \widetilde{X}_T, \Phi \rangle = \int_0^1 \langle X_T(t), \Phi(\cdot, t) \rangle\, dt, \qquad \Phi \in \mathcal{S}(\mathbb{R}^{d+1}).$$

The main ingredients of the proof of the theorem are weak convergence of $\widetilde{X}_T$ and tightness of $\{X_T\}_{T \geq 1}$ in $C([0,1], \mathcal{S}'(\mathbb{R}^d))$.

We need the Laplace transform of $\widetilde{X}_T$. Its form is given in the following lemma.

LEMMA 3.1. *Let $\Phi \in \mathcal{S}(\mathbb{R}^{d+1}), \Phi \geq 0$, and denote*

$$\Psi_T(x,t) = \frac{1}{F_T} \int_{t/T}^1 \Phi(x,s)\, ds.$$

*Then*

$$\text{(3.7)} \qquad Ee^{-\langle \widetilde{X}_T, \Phi \rangle} = \exp\left\{ \int_{\mathbb{R}^d} \int_0^T \Psi_T(x, T-r) v_T(x,r)\, dr\, dx \right. \\ \left. + \frac{V}{1+\beta} \int_{\mathbb{R}^d} \int_0^T v_T^{1+\beta}(x,r)\, dr\, dx \right\},$$

*where $v_T$ satisfies*

$$\text{(3.8)} \qquad v_T(x,t) = \int_0^t \mathcal{T}_{t-r}\left[ \Psi_T(\cdot, T-r)(1 - v_T(\cdot, r)) - \frac{V}{1+\beta} v_T^{1+\beta}(\cdot, r) \right](x)\, dr,$$

$$0 \leq t \leq T.$$

We omit the proof of this lemma because it can be done in the same way as that of (3.23) in [6], using the Feynman–Kac formula and the form of the generating function of the branching law; see (1.1). Note that the function $v_T$ corresponds to $v_{\Psi_T}(x, T - t, t)$ in [6]. By the definition of $v_T$ (see (3.12) and (3.19) in [6]), we have

$$0 \leq v_T \leq 1. \tag{3.9}$$

PROPOSITION 3.2. *Let* $\Phi \in \mathcal{S}(\mathbb{R}^{d+1}), \Phi \geq 0$. *Then*

$$\lim_{T \to \infty} E e^{-\langle \widetilde{X}_T, \Phi \rangle}$$

$$= \exp\left\{ \frac{V}{1+\beta} \int_{\mathbb{R}^d} \int_0^1 \left[ \int_{\mathbb{R}^d} \int_r^1 \Phi(y, s) \right.\right. \tag{3.10}$$

$$\left.\left. \times \int_r^s p_{u-r}(x) \, du \, ds \, dy \right]^{1+\beta} dr \, dx \right\}.$$

PROOF. We assume that $\Phi$ is of the form $\Phi(x, t) = \varphi(x)\psi(t)$, where $\varphi \in \mathcal{S}(\mathbb{R}^d), \psi \in \mathcal{S}(\mathbb{R})$ and $\varphi, \psi \geq 0$. For general $\Phi$, the proof is the same with slightly more complicated notation. Denote

$$\chi(u) = \int_u^1 \psi(s) \, ds, \qquad \chi_T(u) = \chi\left(\frac{u}{T}\right), \qquad \varphi_T(x) = \frac{1}{F_T} \varphi(x). \tag{3.11}$$

Obviously,

$$\chi(u) \leq C \quad \text{and} \quad \chi_T(u) \leq C, \tag{3.12}$$

for some constant $C$.

By Lemma 3.1, the Laplace transform of $\langle \widetilde{X}_T, \Phi \rangle$ can be written as

$$E e^{-\langle \widetilde{X}_T, \Phi \rangle} = \exp\left\{ \frac{V}{1+\beta} I_1(T) + I_2(T) - \frac{V}{1+\beta} I_3(T) \right\}, \tag{3.13}$$

where

$$I_1(T) = \int_{\mathbb{R}^d} \int_0^T \left( \int_0^r \mathcal{T}_{r-u} \varphi_T(x) \chi_T(T-u) \, du \right)^{1+\beta} dr \, dx, \tag{3.14}$$

$$I_2(T) = \int_{\mathbb{R}^d} \int_0^T \varphi_T(x) \chi_T(T-r) v_T(x, r) \, dr \, dx, \tag{3.15}$$

$$I_3(T) = \int_{\mathbb{R}^d} \int_0^T \left[ \left( \int_0^r \mathcal{T}_{r-u} \varphi_T(x) \chi_T(T-u) \, du \right)^{1+\beta} \right.$$
$$\left. - v_T^{1+\beta}(x, r) \right] dr \, dx. \tag{3.16}$$

In (3.15) and (3.16) $v_T$ is the solution of the equation (3.8) with $\Psi_T(x,t) = \varphi_T(x)\chi_T(t)$. Note that by (3.8) and (3.9), $I_1, I_2$ and $I_3$ are nonnegative, the latter fact following from

$$(3.17) \qquad v_T(x,t) \leq \int_0^t \mathcal{T}_{t-u}\varphi_T(x)\chi_T(T-u)\,du.$$

We will show that $I_2(T)$ and $I_3(T)$ converge to 0 as $T \to \infty$, and $\frac{V}{1+\beta}I_1(T)$ converges to the term in the exponent in (3.10).

We consider first $I_1$. Using (3.11), (2.4), and making the changes of variables $u' = \frac{u}{T}$ and $r' = \frac{r}{T}$ in (3.14), we obtain

$$I_1(T) = T^{(d/\alpha)\beta}\int_{\mathbb{R}^d}\int_0^1\left(\int_0^r \mathcal{T}_{T(r-u)}\varphi(x)\chi(1-u)\,du\right)^{1+\beta}dr\,dx$$

$$= T^{(d/\alpha)\beta}\int_{\mathbb{R}^d}\int_0^1\left(\int_0^r\int_{\mathbb{R}^d}p_{T(r-u)}(x-y)\varphi(y)\chi(1-u)\,du\right)^{1+\beta}dr\,dx.$$

Using (3.1) and substituting $x' = xT^{-1/\alpha}$ and $y' = yT^{-1/\alpha}$, we arrive at

$$I_1(T) = \int_{\mathbb{R}^d}\int_0^1\left(\int_0^r\int_{\mathbb{R}^d}p_{r-u}(x-y)\right.$$

$$(3.18) \qquad\qquad \left.\times \chi(1-u)T^{d/\alpha}\varphi(yT^{1/\alpha})\,dy\,du\right)^{1+\beta}dr\,dx$$

$$= \int_0^1\int_{\mathbb{R}^d}(f_r * g_T(x))^{1+\beta}\,dx\,dr,$$

where

$$(3.19)\quad f_r(x) = \int_0^r p_{r-u}(x)\chi(1-u)\,du \quad\text{and}\quad g_T(x) = T^{d/\alpha}\varphi(xT^{1/\alpha}).$$

By (3.12) and (3.3), it follows that $f_r \in L^{1+\beta}(\mathbb{R}^d)$ for any $r \in [0,1]$, therefore, taking into account the form of $g_T$, we have that $f_r * g_T$ converges in $L^{1+\beta}(\mathbb{R}^d)$ to $\int_{\mathbb{R}^d}\varphi(y)\,dy\,f_r$ for any $r \in [0,1]$ as $T \to \infty$. Moreover, by Young's inequality (see, e.g., [28], Theorem 6.3.1), (3.11), (3.19) and (3.3),

$$\sup_{r\in[0,1]}\|f_r * g_T\|_{1+\beta}^{1+\beta} \leq \sup_{r\in[0,1]}\|f_r\|_{1+\beta}^{1+\beta}\|g_T\|_1^{1+\beta}$$

$$\leq C_1\left\|\int_0^1 p_u(\cdot)\,du\right\|_{1+\beta}^{1+\beta}\|\varphi\|_1^{1+\beta} < \infty.$$

We apply the dominated convergence theorem to (3.18) to conclude that

$$(3.20)\quad \lim_{T\to\infty}I_1(T) = \int_0^1\int_{\mathbb{R}^d}\left(\int_{\mathbb{R}^d}\varphi(y)\,dy\right)^{1+\beta}$$

$$\times \left(\int_0^r p_{r-u}(x)\chi(1-u)\,du\right)^{1+\beta}dx\,dr.$$

Recalling the definition of $\chi$ [see (3.11)] and substituting $u' = 1 - u, r' = 1 - r$, it is easy to see that (3.20) is the same as

$$\lim_{T \to \infty} I_1(T) = \int_{\mathbb{R}^d} \int_0^1 \left[ \int_{\mathbb{R}^d} \int_r^1 \varphi(y) \psi(s) \right. \tag{3.21}$$
$$\left. \times \int_r^s p_{u-r}(x) \, du \, ds \, dy \right]^{1+\beta} dr \, dx.$$

Now we proceed to $I_2$.

Applying (3.17) and (3.12) to (3.15), we obtain

$$I_2(T) \leq C_1 \int_{\mathbb{R}^d} \int_0^T \varphi_T(x) \int_0^r \mathcal{T}_{r-u} \varphi_T(x) \, du \, dr \, dx.$$

Substituting $u' = \frac{u}{T}, r' = \frac{r}{T}$ and using (3.11), we obtain

$$I_2(T) \leq C_1 \frac{T^2}{F_T^2} \int_{\mathbb{R}^d} \int_0^1 \int_0^r \varphi(x) \mathcal{T}_{T(r-u)} \varphi(x) \, du \, dr \, dx. \tag{3.22}$$

Next we use (2.4), the Plancherel formula and the fact that $\widehat{\mathcal{T}_s \varphi}(z) = e^{-s|z|^\alpha} \widehat{\varphi}(z)$ ($\widehat{\phantom{x}}$ denoting Fourier transform). Then (3.22) becomes

$$I_2(T) \leq \frac{C_1}{(2\pi)^d} T^{-1+2((d/\alpha)\beta-1)/(1+\beta)}$$
$$\times \int_0^1 \int_{\mathbb{R}^d} \frac{1 - e^{-Tr|x|^\alpha}}{|x|^\alpha} |\widehat{\varphi}(x)|^2 \, dx \, dr. \tag{3.23}$$

It is easy to check that under (2.1) we have $-1 + 2(\frac{d}{\alpha}\beta - 1)/(1+\beta) < 0$, and since also $\alpha < d$, it follows from (3.23) that

$$\lim_{T \to \infty} I_2(T) = 0. \tag{3.24}$$

It remains to prove that $I_3(T)$ also converges to 0.

We need some more notation:

$$J_1(T) = \int_{\mathbb{R}^d} \int_0^T \left[ \int_0^r \mathcal{T}_{r-u}(\varphi_T(\cdot) v_T(\cdot, u))(x) \, du \right]^{1+\beta} dr \, dx, \tag{3.25}$$

$$J_2(T) = \int_{\mathbb{R}^d} \int_0^T \left[ \int_0^r \mathcal{T}_{r-u}(v_T^{1+\beta}(\cdot, u))(x) \, du \right]^{1+\beta} dr \, dx. \tag{3.26}$$

By (3.8) and (3.16),

$$I_3(T) = \int_{\mathbb{R}^d} \int_0^T \left\{ \left[ \int_0^r \mathcal{T}_{r-u} \varphi_T(x) \chi_T(T-u) \, du \right]^{1+\beta} \right.$$
$$- \left[ \int_0^r \mathcal{T}_{r-u} \varphi_T(x) \chi_T(T-u) \, du \right.$$

(3.27)
$$-\int_0^r \mathcal{T}_{r-u}(\varphi_T(\cdot)\chi_T(T-u)v_T(\cdot,u))(x)\,du$$
$$\left.-\frac{V}{1+\beta}\int_0^r \mathcal{T}_{r-u}(v_T^{1+\beta}(\cdot,u))(x)\,du\right]^{1+\beta}\right\}\,dr\,dx.$$

Note that by (3.9) we have

(3.28)
$$\int_0^r \mathcal{T}_{r-u}\varphi_T(x)\chi_T(T-u)\,du$$
$$\geq \int_0^r \mathcal{T}_{r-u}(\varphi_T(\cdot)\chi_T(T-u)v_T(\cdot,u))(x)\,du$$
$$+\frac{V}{1+\beta}\int_0^r \mathcal{T}_{r-u}(v_T^{1+\beta}(\cdot,u))(x)\,du.$$

We apply (3.4) with $\delta = (1+\beta)/2$ and $a = a_1 + a_2$, where
$$a_1 = \int_0^r \mathcal{T}_{r-u}(\varphi_T(\cdot)\chi_T(T-u)v_T(\cdot,u))(x)\,du,$$
$$a_2 = \frac{V}{1+\beta}\int_0^r \mathcal{T}_{r-u}(v_T^{1+\beta}(\cdot,u))(x)\,du,$$
$$a+b = \int_0^r \mathcal{T}_{r-u}\varphi_T(x)\chi_T(T-u)\,du,$$

[$b \geq 0$ by (3.28)], and then we use the estimate $(a_1 + a_2)^{1+\beta} \leq 2^\beta(a_1^{1+\beta} + a_2^{1+\beta})$ to arrive at

$$I_3(T) \leq C_1(J_1(T) + J_2(T))$$
$$+ C_2 \int_{\mathbb{R}^d}\int_0^T \left(\int_0^r \mathcal{T}_{r-u}(\varphi_T(\cdot)\chi_T(T-u)v_T(\cdot,u))(x)\,du\right.$$
$$\left.+\frac{V}{1+\beta}\int_0^r \mathcal{T}_{r-u}(v_T^{1+\beta}(\cdot,u))(x)\,du\right)^{(1+\beta)/2}$$
$$\times \left(\int_0^r \mathcal{T}_{r-u}(\varphi_T(\cdot)\chi_T(T-u)(1-v_T(\cdot,u)))(x)\,du\right.$$
$$\left.-\frac{V}{1+\beta}\int_0^r \mathcal{T}_{r-u}(v_T^{1+\beta}(\cdot,u))(x)\,du\right)^{(1+\beta)/2}\,dr\,dx.$$

Using the estimate $(a+b)^{(1+\beta)/2} \leq a^{(1+\beta)/2} + b^{(1+\beta)/2}$ in the first factor under the integral $\int_{\mathbb{R}^d}\int_0^T$ and (3.28) in the second one, then finally by the Schwarz inequality and (3.12), we obtain

(3.29) $\quad I_3(T) \leq C_1(J_1(T)+J_2(T)) + C_2\sqrt{I_1(T)}(\sqrt{J_1(T)}+\sqrt{J_2(T)}).$

It remains to prove that $J_1(T)$ and $J_2(T)$ tend to 0 as $T \to \infty$.
By (3.17), (3.12) and (3.25), we have

$$J_1(T) \leq C_1 \int_{\mathbb{R}^d} \int_0^T \left[ \int_0^r \mathcal{T}_{r-u}\left( \varphi_T(\cdot) \int_0^u \mathcal{T}_{u-v}\varphi_T(\cdot) \, dv \right)(x) \, du \right]^{1+\beta} dr \, dx$$

$$\leq \frac{T^{3+2\beta}}{F_T^{2+2\beta}} C_1 \int_{\mathbb{R}^d} \left[ \int_0^1 \mathcal{T}_{Tu}\left( \varphi(\cdot) \int_0^1 \mathcal{T}_{Tv}\varphi(\cdot) \, dv \right)(x) \, du \right]^{1+\beta} dx$$

$$= \frac{T^{3+2\beta}}{F_T^{2+2\beta}} C_1 \int_{\mathbb{R}^d} \left[ \int_{\mathbb{R}^{2d}} T^{-2d/\alpha} \int_0^1 p_u((x-y)T^{-1/\alpha})\varphi(y) \right.$$

$$\left. \times \int_0^1 p_v((y-z)T^{-1/\alpha})\varphi(z) \, dv \, du \, dz \, dy \right]^{1+\beta} dx,$$

where we used the definition of $\mathcal{T}_t$ and (3.1) in the last step.

We now recall $F_T$ [see (2.4)] and substitute $x' = xT^{-1/\alpha}, y' = T^{-1/\alpha}y, z' = T^{-1/\alpha}z$ to obtain

$$J_1(T) \leq C_1 T^{-1} \int_{\mathbb{R}^d} \left( \int_{\mathbb{R}^{2d}} \int_0^1 p_u(x-y) \, du \, T^{(d/\alpha)\beta/(1+\beta)} \varphi(T^{1/\alpha}y) \right.$$
(3.30)
$$\left. \times \int_0^1 p_v(y-z) \, dv \, T^{d/\alpha} \varphi(T^{1/\alpha}z) \, dz \, dy \right)^{1+\beta} dx.$$

To simplify the notation, we introduce the functions

$$f(x) = \int_0^1 p_u(x) \, du,$$
(3.31)
$$g_{1,T}(x) = T^{(d/\alpha)\beta/(1+\beta)} \varphi(T^{1/\alpha}x),$$
$$g_{2,T}(x) = T^{d/\alpha} \varphi(T^{1/\alpha}x).$$

It is easy to check that

(3.32) $\|g_{1,T}\|_{(1+\beta)/\beta} = \|\varphi\|_{(1+\beta)/\beta} < \infty$ and $\|g_{2,T}\|_1 = \|\varphi\|_1 < \infty.$

In the notation of (3.31) the inequality (3.30) can be written as

$$J_1(T) \leq C_1 T^{-1} \|f * (g_{1,T}(f * g_{2,T}))\|_{1+\beta}^{1+\beta}.$$

We use consecutively the Young, the Hölder and again the Young inequalities, obtaining

$$J_1(T) \leq C_1 T^{-1} \|f\|_{1+\beta}^{1+\beta} \|g_{1,T}(f * g_{2,T})\|_1^{1+\beta}$$
$$\leq C_1 T^{-1} \|f\|_{1+\beta}^{1+\beta} \|g_{1,T}\|_{(1+\beta)/\beta}^{1+\beta} \|f * g_{2,T}\|_{1+\beta}^{1+\beta}$$
$$\leq C_1 T^{-1} \|f\|_{1+\beta}^{1+\beta} \|g_{1,T}\|_{(1+\beta)/\beta}^{1+\beta} \|f\|_{1+\beta}^{1+\beta} \|g_{2,T}\|_1^{1+\beta}.$$

By (3.31), (3.32) and (3.3), it follows that

(3.33) $$\lim_{T\to\infty} J_1(T) = 0.$$

The term $J_2(T)$ can be dealt with in a similar manner. By (3.26), (3.17) and (3.12), we have

$$J_2(T) \leq C_1 \int_{\mathbb{R}^d} \int_0^T dr \left( \int_0^r \mathcal{T}_{r-u} \left( \int_0^u \mathcal{T}_{u-v} \varphi_T \, dv \right)^{1+\beta} (x) \, du \right)^{1+\beta} dx$$

$$\leq C_1 \frac{T^{1+(1+\beta)+(1+\beta)^2}}{F_T^{(1+\beta)(1+\beta)}} \int_{\mathbb{R}^d} \left( \int_0^1 \mathcal{T}_{Tu} \left( \int_0^1 \mathcal{T}_{Tv} \varphi \, dv \right)^{1+\beta} (x) \, du \right)^{1+\beta} dx.$$

As in the case of $J_1$, we use consecutively the definition of $\mathcal{T}_t$, (3.1), substitutions $x' = xT^{-1/\alpha}, y' = yT^{-1/\alpha}, z' = zT^{-1/\alpha}$, and (2.4) to obtain

$$J_2(T) \leq C_1 T^{1-d/\alpha\beta}$$

(3.34)
$$\times \int_{\mathbb{R}^d} \left[ \int_{\mathbb{R}^d} \int_0^1 p_u(x-y) \, du \right.$$
$$\left. \times \left( \int_{\mathbb{R}^d} \int_0^1 p_v(y-z) \, dv T^{d/\alpha} \varphi(T^{1/\alpha} z) \, dz \right)^{1+\beta} dy \right]^{1+\beta} dx.$$

Using notation (3.31) and applying the Young inequality twice, (3.34) can be estimated in as follows:

$$J_2(T) \leq C_1 T^{1-(d/\alpha)\beta} \| f * (f * g_{2,T})^{1+\beta} \|_{1+\beta}^{1+\beta}$$

$$\leq C_1 T^{1-(d/\alpha)\beta} \| f \|_{1+\beta}^{1+\beta} \| (f * g_{2,T})^{1+\beta} \|_1^{1+\beta}$$

$$\leq C_1 T^{1-(d/\alpha)\beta} \| f \|_{1+\beta}^{1+\beta} \| f * g_{2,T} \|_{1+\beta}^{(1+\beta)(1+\beta)}$$

$$\leq C_1 T^{1-(d/\alpha)\beta} \| f \|_{1+\beta}^{1+\beta} \| f \|_{1+\beta}^{(1+\beta)(1+\beta)} \| g_{2,T} \|_1^{(1+\beta)(1+\beta)}.$$

By (3.32) and (3.3), we obtain

(3.35) $$\lim_{T\to\infty} J_2(T) = 0.$$

By (3.33), (3.35), (3.24) and (3.29), we get

(3.36) $$\lim_{T\to\infty} I_3(T) = 0.$$

Putting together (3.13), (3.21), (3.24) and (3.36) finishes the proof of the proposition. □

We now pass to the tightness. We state a slightly more general result which includes also the lower critical dimension $d = \alpha/\beta$. This is used in a forthcoming paper [30].

LONG RANGE DEPENDENCE AND BRANCHING    17PROPOSITION 3.3. *Assume that*

(3.37) $$\frac{\alpha}{\beta} \leq d < \frac{\alpha(1+\beta)}{\beta}.$$

*Then the family $\{X_T\}_{T \geq 1}$ is tight in $C([0,1], \mathcal{S}'(\mathbb{R}^d))$.*

PROOF. The fact that the process $X_T$ lacks moments of order $\geq 1 + \beta$ for $\beta < 1$ prevents the use of standard methods for proving tightness. Also, the Laplace transform technique we have employed for showing weak convergence of the space-time random field $\widetilde{X}_T$ does not seem to be amenable to a tightness proof in the present case. Instead, we will give a proof based on the characteristic function of $\widetilde{X}_T$.

By Theorem 12.3 of Billingsely [1] and the theorem of Mitoma [24], it suffices to show that for any $\varphi \in \mathcal{S}(\mathbb{R}^d)$ there exist constants $\nu \geq 0$ and $\gamma > 0$ such that

(3.38) $$P(|\langle X_T(t_2), \varphi \rangle - \langle X_T(t_1), \varphi \rangle| \geq \delta) \leq \frac{C(\varphi)}{\delta^\nu}(t_2 - t_1)^{1+\gamma}$$

holds for all $t_1, t_2 \in [0,1]$, $t_1 < t_2$, all $T \geq 1$, and all $\delta > 0$.

Since each $\varphi \in \mathcal{S}(\mathbb{R}^d)$ can be written as $\varphi = \varphi_1 - \varphi_2$, $\varphi_1, \varphi_2 \in \mathcal{S}(\mathbb{R}^d)$, $\varphi_1, \varphi_2 \geq 0$, (see, e.g., the lemma in Section 3 of [6]), it suffices to assume $\varphi \geq 0$, which we do from now on. So, fix $\varphi \geq 0$ and $t_1, t_2$.

In order to prove (3.38), we use the estimate (see, e.g., [2], Proposition 8.29)

(3.39) $$P(|\langle \widetilde{X}_T, \varphi \otimes \psi \rangle| \geq \delta)$$
$$\leq C\delta \int_0^{1/\delta} (1 - \text{Re}(E \exp\{-i\theta \langle \widetilde{X}_T, \varphi \otimes \psi \rangle\})) \, d\theta,$$

where (arguing as in the tightness proof in [7]) $\psi \in \mathcal{S}(\mathbb{R})$ is an approximation of $\delta_{t_2} - \delta_{t_1}$ supported on $[t_1, t_2]$ such that

$$\chi(t) = \int_t^1 \psi(s) \, ds$$

satisfies

(3.40) $$\chi \in \mathcal{S}(\mathbb{R}), \qquad 0 \leq \chi \leq \mathbb{1}_{[t_1, t_2]}.$$

Hence, it suffices to prove that for any $\chi$ satisfying (3.40), the right-hand side of (3.39) (with the corresponding $\psi$) is estimated from above by the right-hand side of (3.38), with constants not depending on $\chi$. To this end, we define a complex-valued analogue of the function $v_T$ considered before, namely (using the same notation),

$$(3.41) \quad v_{\theta,T}(x,t) = 1 - E\exp\left\{-i\theta \int_0^t \langle N_s^x, \varphi_T\rangle \chi_T(T-t+s)\,ds\right\}, \qquad \theta > 0,$$

where $\varphi_T$, $\chi_T$ are given by (3.11), and $N^x$ is the empirical measure of the branching system started from a single particle at $x$ (see (3.12) and (3.19) in [6]).

Since the Feyman–Kac formula holds for complex-valued functions, the same procedure used before (see Lemma 3.1) shows that $v_{\theta,T}$ satisfies

$$(3.42) \quad v_{\theta,T}(x,t) = \int_0^t \mathcal{T}_{t-s}\left[i\theta \varphi_T \chi_T(T-s)(1-v_{\theta,T}(\cdot,s)) - \frac{V}{1+\beta}v_{\theta,T}^{1+\beta}(\cdot,s)\right](x)\,ds$$

[cf. (3.8)], and

$$(3.43) \quad E\exp\{-i\theta\langle \widetilde{X}_T, \varphi\otimes\psi\rangle\} = \exp\{I+II\}, \qquad \theta > 0,$$

where

$$(3.44) \quad I = i\theta \int_{\mathbb{R}^d}\int_0^T \varphi_T(x)\chi_T(T-s)v_{\theta,T}(x,s)\,ds\,dx,$$

$$(3.45) \quad II = \frac{V}{1+\beta}\int_{\mathbb{R}^d}\int_0^T v_{\theta,T}^{1+\beta}(x,s)\,ds\,dx$$

[cf. (3.7)]. In these equations $z^{1+\beta} = \exp\{(1+\beta)\log z\}$ is understood in the sense of the principal branch of the logarithm.

Using the inequality

$$(3.46) \quad |1-e^z| \le 2|z|$$

if $|e^z| \le 1$, $z \in \mathbb{C}$, we have, by (3.43),

$$(3.47) \quad 0 \le 1 - \mathrm{Re}\,E\exp\{-i\theta\langle \widetilde{X}_T,\varphi\otimes\psi\rangle\} \le 2(|I|+|II|).$$

Now, by virtue of the previous discussion, taking into account (3.39) and (3.47), we see that (3.38) will be proved if we show that

$$(3.48) \quad |I| \le C(\varphi,\gamma)\theta^2(t_2-t_1)^{1+\gamma}$$

and

$$(3.49) \quad |II| \le C(\varphi,\gamma,V,\beta)\theta^{1+\beta}(t_2-t_1)^{1+\gamma}$$

for any $\gamma$ such that

(3.50) $$0 < \gamma < 1 + \beta - \frac{d\beta}{\alpha}.$$

To prove (3.48) and (3.49), first we observe that by (3.41) and (3.46) we have

(3.51) $$|v_{\theta,T}(x,t)| \leq 2\theta E \int_0^t \langle N_s^x, \varphi_T \rangle \chi_T(T - t + s)\, ds$$
$$= 2\theta \int_0^t \mathcal{T}_{t-s}\varphi_T(x)\chi_T(T - s)\, ds,$$

since $E\langle N_s^x, \varphi \rangle = \mathcal{T}_s\varphi(x)$ (this known fact is obtained by the usual renewal argument).

Combining this and (3.44), we obtain

(3.52) $$|I| \leq \frac{2\theta^2}{F_T^2} \int_0^T \int_0^s \int_{\mathbb{R}^d} \varphi(x) \mathcal{T}_{s-r}\varphi(x)\, dx\, \chi\left(1 - \frac{r}{T}\right)\chi\left(1 - \frac{s}{T}\right) dr\, ds$$
$$= \frac{2\theta^2}{(2\pi)^d} \frac{T^2}{F_T^2} \int_0^1 \chi(s) \int_{\mathbb{R}^d} |\widehat{\varphi}(z)|^2 \int_s^1 e^{-T(r-s)|z|^\alpha} \chi(r)\, dr\, dz\, ds,$$

where, besides obvious substitutions, we have used the Plancherel formula as in (3.22).

Fix $\gamma$ satisfying (3.50) and note that $\frac{1}{\gamma} > 1$ by (3.37).

By Hölder's inequality, we have

(3.53) $$\int_s^1 e^{-T(r-s)|z|^\alpha} \chi(r)\, dr$$
$$\leq \left(\int_s^1 e^{-(T/(1-\gamma))(r-s)|z|^\alpha}\, dr\right)^{1-\gamma} \left(\int_s^1 \chi^{1/\gamma}(r)\, dr\right)^\gamma$$
$$\leq (1 - \gamma)^{1-\gamma} T^{-1+\gamma} |z|^{-\alpha(1-\gamma)} (t_2 - t_1)^\gamma,$$

where in the last estimate we used (3.40).

Combining (3.52) and (3.53) and using (3.40), once again we obtain

$$|I| \leq C(\gamma)\theta^2 \frac{T^{1+\gamma}}{F_T^2} \int_{\mathbb{R}^d} |\widehat{\varphi}(z)|^2 |z|^{-\alpha(1-\gamma)}\, dz (t_2 - t_1)^{1+\gamma}$$
$$= C(\varphi, \gamma)\theta^2 \frac{T^{1+\gamma}}{F_T^2} (t_2 - t_1)^{1+\gamma},$$

and this will imply (3.48) if we show that $T^{1+\gamma}/F_T^2$ is bounded in $T > 1$. To this end, by (2.4), we need to check that

$$1 + \gamma - \frac{2(2 + \beta - (d\beta)/\alpha)}{1 + \beta} \leq 0,$$

or
$$\gamma \leq \frac{3+\beta-2\beta^{d/\alpha}}{1+\beta},$$

and this inequality is indeed satisfied by (3.50) since it can be easily verified [using (3.37)] that

$$1+\beta-\frac{d}{\alpha}\beta \leq \frac{3+\beta-2\beta\frac{d}{\alpha}}{1+\beta}.$$

We now pass to the proof of (3.49). By (3.45) and (3.51), we have

$$|II| \leq \frac{V}{1+\beta}2^{1+\beta}\theta^{1+\beta}\int_{\mathbb{R}^d}\int_0^T\left(\int_0^s \mathcal{T}_{s-r}\varphi_T(x)\chi_T(T-r)\,dr\right)^{1+\beta}ds\,dx$$

[by (3.1), (2.4) and (3.11)]

$$= C(V,\beta)\theta^{1+\beta}\int_{\mathbb{R}^d}\int_0^1 T^{d/\alpha\beta}\left(\int_0^s\int_{\mathbb{R}^d}T^{-d/\alpha}p_{s-r}(T^{-1/\alpha}x-T^{-1/\alpha}y)\right.$$

$$\left.\times \varphi(y)\chi(1-r)\,dy\,dr\right)^{1+\beta}ds\,dx$$

$$= C(V,\beta)\theta^{1+\beta}\int_0^1\int_{\mathbb{R}^d}\left(\int_0^s\int_{\mathbb{R}^d}p_{s-r}(x-y)T^{d/\alpha}\varphi(T^{1/\alpha}y)\right.$$

$$\left.\times \chi(1-r)\,dy\,dr\right)^{1+\beta}ds\,dx$$

$$\leq C(V,\beta)\theta^{1+\beta}\int_0^1\left\|\int_0^s p_{s-r}(\cdot)\chi(1-r)\,dr\right\|_{1+\beta}^{1+\beta}ds\|T^{d/\alpha}\varphi(T^{1/\alpha}\cdot)\|_1^{1+\beta},$$

where we have used the Young inequality in the last step. It is now clear that (3.49) will be proved if we show that

(3.54) $$\int_{\mathbb{R}^d}\int_0^1\left(\int_0^s p_{s-r}(x)\chi(1-r)\,dr\right)^{1+\beta}ds\,dx \leq C(t_2-t_1)^{1+\gamma}.$$

Denote $f_x(s) = p_s(x)\mathbb{1}_{[0,1]}(s)$ for any $x \in \mathbb{R}^d$, and $g(s) = \chi(1-s)$, and put

(3.55) $$p = \frac{1+\beta}{1+\beta-\gamma}, \qquad q = \frac{1+\beta}{1+\gamma}.$$

Note that by (3.50),

(3.56) $$1 < p < \frac{\alpha}{d\beta}(1+\beta), \qquad 1 < q,$$

and

$$\frac{1}{p}+\frac{1}{q}-\frac{1}{1+\beta}=1.$$

By the Young inequality, we have

$$\int_0^1 \left( \int_0^s p_{s-r}(x) \chi(1-r)\, dr \right)^{1+\beta} ds \leq \|f_x * g\|_{1+\beta}^{1+\beta} \leq \|f_x\|_p^{1+\beta} \|g\|_q^{1+\beta}.$$

But $\|g\|_q^{1+\beta} \leq (t_2 - t_1)^{1+\gamma}$ by (3.40) and (3.55), therefore, to finish the proof of (3.54) (and the proposition as well), it suffices to show that

$$\int_{\mathbb{R}^d} \|f_x\|_p^{1+\beta}\, dx = A + B < \infty,$$

where

$$A = \int_{|x|>1} \|f_x\|_p^{1+\beta}\, dx, \qquad B = \int_{|x|\leq 1} \|f_x\|_p^{1+\beta}\, dx.$$

By (3.1) and (3.2), we have

$$A = \int_{|x|>1} \left( \int_0^1 s^{-pd/\alpha} p_1^p(s^{-1/\alpha} x)\, ds \right)^{(1+\beta)/p} dx$$

$$\leq C \int_{|x|>1} |x|^{-(1+\beta)(d+\alpha)}\, dx \left( \int_0^1 s^{-pd/\alpha + p(1+d/\alpha)}\, ds \right)^{(1+\beta)/p} < \infty.$$

$B$ is estimated analogously. We observe that from (3.56) we have

$$0 < \frac{(d/\alpha)p - 1}{((d+\alpha)/\alpha)p} < \frac{d}{(d+\alpha)(1+\beta)}.$$

So, if we fix $\rho$ such that

$$\frac{(d/\alpha)p - 1}{((d+\alpha)/\alpha)p} < \rho < \frac{d}{(d+\alpha)(1+\beta)},$$

and use the fact that $p_1(x) \leq C(\rho) p_1^\rho(x)$ (since $p_1$ is bounded and $\rho \leq 1$), then we have, by (3.1) and (3.2),

$$B \leq C \int_{|x|\leq 1} |x|^{-(1+\beta)(d+\alpha)\rho}\, dx \left( \int_0^1 s^{-pd/\alpha + p(1+d/\alpha)\rho}\, ds \right)^{(1+\beta)/p} < \infty. \quad \square$$

We still need a result on convergence of Laplace transforms, which we formulate as a lemma. This result is known, following from a standard argument on analytic extensions (see, e.g., [20] for the proof in the case $d = 1$, the proof for $d = 2$ is analogous).

LEMMA 3.4. *Let $0 < \beta \leq 1$.*

(a) *If $\eta_n, n = 1, 2, \ldots$ are real random variables such that*

(3.57) $$\lim_{n \to \infty} E e^{-\rho \eta_n} = e^{K \rho^{1+\beta}}, \qquad \rho > 0,$$

*for some constant $K > 0$, then $\eta_n$ converges in distribution to a random variable $\eta$ whose law is $(1+\beta)$-stable, totally skewed to the right, with characteristic function*

$$Ee^{iz\eta} = \exp\left\{K\left(\cos\frac{\pi}{2}(1+\beta)\right)|z|^{1+\beta}\left(1 - i(\operatorname{sgn} z)\tan\frac{\pi}{2}(1+\beta)\right)\right\},$$
(3.58)
$$z \in \mathbb{R}.$$

(b) *If $\eta_n, n = 1, 2, \ldots$ are two-dimensional random variables and $\mu$ is a finite measure on $S_+^2 = \{(u_1, u_2) : u_1^2 + u_2^2 = 1, u_1, u_2 \geq 0\}$ such that*

$$(3.59) \qquad \lim_{n\to\infty} Ee^{-\rho\cdot\eta_n} = \exp\left\{\int_{S_+^2}(\rho\cdot y)^{1+\beta}\mu(dy)\right\}, \qquad \rho \in \mathbb{R}_+^2,$$

*then $\eta_n$ converges in distribution to a $(1+\beta)$-stable random variable $\eta$ with characteristic function*

$$Ee^{iz\cdot\eta} = \exp\left\{\left(\cos\frac{\pi}{2}(1+\beta)\right)\int_{S_+^2}|z\cdot y|^{1+\beta}\right.$$
(3.60)
$$\left.\times\left(1 - i(\operatorname{sgn}(z\cdot y))\tan\frac{\pi}{2}(1+\beta)\right)\mu(dy)\right\},$$
$$z \in \mathbb{R}^2$$

*($\cdot$ denotes the inner product in $\mathbb{R}^2$).*

COROLLARY 3.5. *For each $\Phi \in \mathcal{S}(\mathbb{R}^{d+1})$,*

$$(3.61) \qquad \langle \widetilde{X}_T, \Phi \rangle \Rightarrow \xi_\Phi \qquad \text{as } T \to \infty,$$

*where $\xi_\Phi$ is $(1+\beta)$-stable with*

$$Ee^{i\xi_\Phi} = \exp\left\{-K_1\int_{\mathbb{R}^d}\int_0^1\left|\int_{\mathbb{R}^d}\int_r^1\Phi(y,s)\int_r^s p_{u-r}(x)\,du\,ds\,dy\right|^{1+\beta}\right.$$
(3.62)
$$\times\left(1 - i\operatorname{sgn}\left(\int_{\mathbb{R}^d}\int_r^1\Phi(y,s)\int_r^s p_{u-r}(x)\,du\,ds\,dy\right)\right.$$
$$\left.\left.\times\tan\frac{\pi}{2}(1+\beta)\right)dr\,dx\right\},$$

*where*

$$(3.63) \qquad K_1 = -\frac{V}{1+\beta}\cos\frac{\pi}{2}(1+\beta).$$

PROOF. Proposition 3.2 and Lemma 3.4(a) imply that (3.61) holds for any $\Phi \geq 0$. This, applied to $\rho_1 \Phi_1 + \rho_2 \Phi_2$, $\rho_1, \rho_2 \geq 0, \Phi_1, \Phi_2 \geq 0$, and Lemma 3.4(b) imply weak convergence of $(\langle \widetilde{X}_T, \Phi_1 \rangle, \langle \widetilde{X}_T, \Phi_2 \rangle)$ with a $\mu$ obtained in a standard way (see, e.g., [20], proof of Theorem 5.6). Hence, (3.61) follows since, analogously as in the proof of Proposition 3.3, we can write any $\Phi \in \mathcal{S}(\mathbb{R}^{d+1})$ as $\Phi = \Phi_1 - \Phi_2, \Phi_1, \Phi_2 \geq 0$. $\square$

Now we can complete the proof of Theorem 2.2. By Proposition 3.3, we know that $\{X_T\}_{T \geq 1}$ is tight. Let $X_{T_n} \Rightarrow X$ in $C([0,1], \mathcal{S}'(\mathbb{R}^d))$ for some $T_n \nearrow \infty$. Then $\langle \widetilde{X}_{T_n}, \Phi \rangle \Rightarrow \langle \widetilde{X}, \Phi \rangle$ for any $\Phi \in \mathcal{S}(\mathbb{R}^{d+1})$ (see [4]), and

$$E \exp\{i\langle \widetilde{X}, \Phi \rangle\} = E \exp\{i\xi_\Phi\}.$$

We now find the finite dimensional distributions of $X$ arguing, for instance, as in the proof of Proposition 4.1 of [4]. Fix $0 \leq t_1 < \cdots < t_k \leq 1$, $z_1, \ldots, z_k \in \mathbb{R}$, $\varphi_1, \ldots, \varphi_k \in \mathcal{S}(\mathbb{R}^d)$.

Let $\psi_{m,j} \in \mathcal{S}(\mathbb{R})$, $\psi_{m,j} \to \delta_{t_j}$ in $\mathcal{S}'(\mathbb{R})$ as $m \to \infty$, $j = 1, \ldots, k$. Define

$$\Phi_m(x,t) = \sum_{j=1}^{k} z_j \varphi_j(x) \psi_{m,j}(t);$$

then

$$\lim_{m \to \infty} \langle \widetilde{X}, \Phi_m \rangle = \sum_{j=1}^{k} z_j \langle X(t_j), \varphi_j \rangle,$$

and it is easily seen that the right-hand side of (3.62) (with $\Phi_m$ instead of $\Phi$) converges as $m \to \infty$ to

$$E \exp\{i K_1^{1/(1+\beta)} (z_1 \langle \lambda, \varphi_1 \rangle \xi_{t_1} + \cdots + z_k \langle \lambda, \varphi_k \rangle \xi_{t_k})\}$$

[see (3.63)], where the distributions of $\xi$ are given by (2.5). This means that $X = K_1^{1/(1+\beta)} \lambda \xi$, and Theorem 2.2 is proved.

REMARK 3.6. The proof of Theorem 2.2 is also valid for $\beta = 1$, so we obtain the result of [6] as a special case. Note, however, that, in contrast to [6], in this proof the form of the covariance of the empirical process $N_t$ is not needed.

**4. Proof of Theorem 2.7.** By (2.5), it is clear that it suffices to investigate the asymptotics of

(4.1) $$D_T^+ = D_T(1, z; u, v, s, t), \qquad z > 0,$$

and

(4.2) $$D_T^- = D_T(1, -z; u, v, s, t), \qquad z > 0.$$

Fix $0 \leq u < v < s < t$ and $z > 0$. The theorem will be proved if we show that

(4.3) $\quad D_T^\pm \leq CT^{-d/\alpha} \quad$ if either $\alpha = 2$, or $\alpha < 2$ and $\beta > \dfrac{d}{d+\alpha}$,

(4.4)
$$D_T^\pm \leq CT^{-(d/\alpha)\delta} \quad \text{for any } \beta < \delta < 1 + \beta - \dfrac{d}{d+\alpha}$$
$$\text{if } \alpha < 2,\ \beta \leq \dfrac{d}{d+\alpha},$$

and for $T$ sufficiently large,

(4.5) $\quad D_T^+ \geq CT^{-d/\alpha}$,

(4.6) $\quad D_T^+ \geq CT^{-(d/\alpha)\delta} \quad$ for any $\delta > 1 + \beta - \dfrac{d}{d+\alpha}$ if $\alpha < 2$, $\beta \leq \dfrac{d}{d+\alpha}$.

Here and in the sequel the constants $C, C_1$, and so forth are different in each line and may depend on $d, \alpha, \beta, u, v, s, t, z$, but never on $T$.

Denote

(4.7)
$$\begin{aligned}
U &= z\mathbb{1}_{[0,t+T]}(r) \int_r^{t+T} p_{r'-r}(x)\, dr' \\
&\quad - z\mathbb{1}_{[0,s+T]}(r) \int_r^{s+T} p_{r'-r}(x)\, dr' \\
&= z\Bigg( \mathbb{1}_{[0,s+T]}(r) \int_{s+T}^{t+T} p_{r'-r}(x)\, dr' \\
&\quad + \mathbb{1}_{(s+T,t+T]}(r) \int_r^{t+T} p_{r'-r}(x)\, dr' \Bigg),
\end{aligned}$$

(4.8)
$$\begin{aligned}
R &= \mathbb{1}_{[0,v]}(r) \int_r^v p_{r'-r}(x)\, dr' - \mathbb{1}_{[0,u]}(r) \int_r^u p_{r'-r}(x)\, dr' \\
&= \mathbb{1}_{[0,u]}(r) \int_u^v p_{r'-r}(x)\, dr' + \mathbb{1}_{(u,v]}(r) \int_r^v p_{r'-r}(x)\, dr'.
\end{aligned}$$

By (2.5), (2.6) and (4.1), we have

(4.9)
$$D_T^+ = \left| \left(1 - i\tan\dfrac{\pi}{2}(1+\beta)\right) \int_{\mathbb{R}^{d+1}} [(U+R)^{1+\beta} - U^{1+\beta} - R^{1+\beta}]\, dx\, dr \right|,$$

and, analogously,

(4.10)
$$\begin{aligned}
D_T^- &= \bigg| \int_{\mathbb{R}^{d+1}} \Big[ |R-U|^{1+\beta} - U^{1+\beta} - R^{1+\beta} \Big] dx\, dr \\
&\quad + i\left(\tan\dfrac{\pi}{2}(1+\beta)\right) \int_{\mathbb{R}^{d+1}} [-|R-U|^{1+\beta} \operatorname{sgn}(R-U)
\end{aligned}$$

$$- U^{1+\beta} + R^{1+\beta}] \, dx \, dr \bigg|.$$

Denote for brevity

$$(4.11) \qquad f = f(x,r) = z \int_{s+T}^{t+T} p_{r'-r}(x) \, dr',$$

$$(4.12) \qquad g_1 = g_1(x,r) = \int_u^v p_{r'-r}(x) \, dr',$$

$$(4.13) \qquad g_2 = g_2(x,r) = \int_r^v p_{r'-r}(x) \, dr'.$$

By (3.1) and (3.2), we have

$$(4.14) \quad f(x,r) = z \int_{s-r}^{t-r} (r'+T)^{-d/\alpha} p_1((r'+T)^{-1/\alpha} x) \, dr' \leq CT^{-d/\alpha},$$

$$(4.15) \quad g_j(x,r) \leq \frac{C}{|x|^{d+\alpha}}, \qquad j = 1, 2,$$

$$(4.16) \quad g_j(x,r) \leq \int_0^v p_{r'}(x) \, dr', \qquad j = 1, 2.$$

Note that the constants $C$ in (4.14) and (4.15) do not depend on $x, r$.

By (4.7)–(4.9),

$$(4.17) \quad \begin{aligned} D_T^+ = C \bigg[ &\int_0^u \int_{\mathbb{R}^d} ((f+g_1)^{1+\beta} - f^{1+\beta} - g_1^{1+\beta}) \, dx \, dr \\ &+ \int_u^v \int_{\mathbb{R}^d} ((f+g_2)^{1+\beta} - f^{1+\beta} - g_2^{1+\beta}) \, dx \, dr. \bigg] \end{aligned}$$

Using (3.4) with $\delta = 1$, we obtain (omitting $dx \, dr$)

$$D_T^+ \leq C(1+\beta) \bigg[ \int_0^u \int_{|x| \leq 1} f g_1^\beta + \int_0^u \int_{|x| > 1} f g_1^\beta \\ + \int_u^v \int_{|x| \leq 1} f g_2^\beta + \int_u^v \int_{|x| > 1} f g_2^\beta \bigg].$$

Assume $\alpha < 2$ and $\beta > d/(d+\alpha)$. Then (4.3) for $D_T^+$ follows by (4.14), (4.15) (for $|x| > 1$) and (4.16) (for $|x| \leq 1$). If $\alpha = 2$, then $g_j^\beta$ is integrable over $\{|x| > 1\}$ for all $\beta$ since by the properties of the Brownian density, we have $g_j(x,r) \leq C_1 e^{-C_2|x|}$ for $|x| \geq 1$ [instead of (4.15)], so we obtain (4.3) without any restriction on $\beta$.

Now assume $\alpha < 2$ and $\beta \leq d/(d+\alpha)$. Fix $\beta < \delta < 1 + \beta - d/(d+\alpha)$, and apply (3.4) to (4.17). Again by (4.14)–(4.16), we have

$$D_T^+ \leq C_1 T^{-(d/\alpha)\delta} \bigg[ u \int_{|x| \leq 1} \bigg( \int_0^v p_{r'}(x) \, dr' \bigg)^{1+\beta-\delta} dx$$

$$+ u \int_{|x|>1} \frac{1}{|x|^{(d+\alpha)(1+\beta-\delta)}}\, dx$$

$$+ (v-u) \int_{|x|\leq 1} \left( \int_0^v p_{r'}(x)\, dr' \right)^{1+\beta-\delta} dx$$

$$+ (v-u) \int_{|x|>1} \frac{1}{|x|^{(d+\alpha)(1+\beta-\delta)}}\, dx \bigg].$$

Hence, (4.4) for $D_T^+$ follows since $1+\beta-\delta < 1$ and $(d+\alpha)(1+\beta-\delta) > d$.

(4.3) and (4.4) for $D_T^-$ are derived in the same way. We only apply the following easy consequences of (3.4):

$$||a-b|^{1+\beta} - a^{1+\beta} - b^{1+\beta}| \leq (3+\beta) a^\delta b^{1+\beta-\delta},$$
$$\beta \leq \delta \leq 1,\ a,b \geq 0,$$

$$||a-b|^{1+\beta}\mathrm{sgn}(a-b) + b^{1+\beta} - a^{1+\beta}| \leq (1+\beta) a^\delta b^{1+\beta-\delta},$$
$$\beta \leq \delta \leq 1,\ a,b, \geq 0.$$

We now pass to the lower estimates of $D_T^+$.

By (4.17),

$$(4.18) \quad D_T^+ \geq C \int_u^{(u+v)/2} \int_{|x|\leq 1} ((f+g_2)^{1+\beta} - f^{1+\beta} - g_2^{1+\beta})\, dx\, dr.$$

Note that for $g_2$ defined by (4.13), if $|x| \leq 1$ and $r \in [u, \frac{u+v}{2}]$, we have

$$g_2(x,r) \geq \int_{(v-u)/4}^{(v-u)/2} p_{r'}(x)\, dr' \geq C(u,v,d,\alpha) > 0,$$

by (3.1). This, combined with (4.14), implies

(4.19) $\qquad\qquad\qquad g_2 \geq f \qquad$ for $T$ large,

for $x,r$ as above. Hence, we can apply (3.5) to (4.18) and obtain

$$D_T^+ \geq C \int_u^{(u+v)/2} \int_{|x|\leq 1} f(x,r) g_2^\beta(x,r)\, dx\, dr$$
(4.20)
$$\geq C_1 \int_u^{(u+v)/2} \int_{|x|\leq 1} f(x,r)\, dx\, dr.$$

For $x,r$ as above, using the equality in (4.14), we have

$$f(x,r) \geq z \int_{s-(u+v)/2}^{t-u} (r'+T)^{-d/\alpha} p_1((r'+T)^{-1/\alpha} x)\, dr'$$
(4.21)
$$\geq C T^{-d/\alpha},$$

thus obtaining (4.5).

Finally, assume $\alpha < 2$ and $\beta \leq d/(d+\alpha)$. Let $0 < \varepsilon < d/(d+\alpha)\alpha$.

By (4.17),

$$(4.22) \quad D_T^+ \geq C \int_u^{(u+v)/2} \int_{1 \leq |x| \leq T^{d/((d+\alpha)\alpha)-\varepsilon}} ((f+g_2)^{1+\beta} - f^{1+\beta} - g_2^{1+\beta}) \, dx \, dr.$$

For $x, r$ in the domain of integration in (4.22), we have

$$(4.23) \quad \begin{aligned} g_2(x,r) &\geq \int_{(v-u)/4}^{(v-u)/2} (r')^{-d/\alpha} p_1((r')^{-1/\alpha} x) \, dr' \\ &\geq \left(\frac{v-u}{2}\right)^{-d/\alpha} \frac{v-u}{4} p_1\left(\left(\frac{v-u}{4}\right)^{-1/\alpha} x\right) \\ &\geq C \frac{1}{|x|^{d+\alpha}} \\ &\geq C T^{-d/\alpha + \varepsilon(d+\alpha)} \end{aligned}$$

by (3.1) and (3.2). Taking into account the equality in (4.14) and (3.2), we see that (4.19) holds also in the present case. Inequality (3.5) applied to (4.22) yields

$$(4.24) \quad D_T^+ \geq C_1 \int_u^{(u+v)/2} \int_{1 \leq |x| \leq T^{d/((d+\alpha)\alpha)-\varepsilon}} f(x,r) g_2^\beta(x,r) \, dx \, dr.$$

As $p_1(T^{-1/\alpha} x) \geq C_2 > 0$ for $T > 1$ and $|x| \leq T^{d/((d+\alpha)\alpha)-\varepsilon}$, from the equality in (4.14), we obtain (4.21) again.

Combining the estimates (4.23), (4.24) and (4.21), we see that, for large $T$,

$$\begin{aligned} D_T^+ &\geq C T^{(-d/\alpha + \varepsilon(d+\alpha))\beta} T^{-d/\alpha} T^{(d/((d+\alpha)\alpha)-\varepsilon)d} \\ &= C T^{-(d/\alpha)[\beta + 1 - d/(d+\alpha) + \varepsilon\alpha(1-\beta(d+\alpha)/d)]}. \end{aligned}$$

This implies (4.6) since $1 - \beta \frac{d+\alpha}{d} \geq 0$ and $\varepsilon$ can be made arbitrarily small.

**Acknowledgments.** We thank a referee for a careful reading of the paper and comments which led to improvements. We also thank the hospitality of the Institute of Mathematics, National University of Mexico (UNAM), where this paper was completed.

LONG RANGE DEPENDENCE AND BRANCHING 29[22] MARCUS, M. B. and ROSIŃSKI, J. (2005). Continuity and boundedness of infinitely divisible processes: A Poisson point process approach. *J. Theoret. Probab.* **18** 109–160. MR2132274

[23] MIŁOŚ, P. (2007). Occupation time fluctuations of Poisson and equilibrium finite variance branching systems. *Probab. Math. Statist.* To appear.

[24] MITOMA, I. (1983). Tightness of probabilities on $C([0,1]; \mathcal{S}')$ and $D([0,1]; \mathcal{S}')$. *Ann. Probab.* **11** 989–999. MR0714961

[25] ROSIŃSKI, J. and ŻAK, T. (1997). The equivalence of ergodicity and weak mixing for infinitely divisible processes. *J. Theoret. Probab.* **10** 73–86. MR1432616

[26] SAMORODNITSKY, G. and TAQQU, M. S. (1994). *Stable Non-Gaussian Random Processes*. Chapman and Hall, New York. MR1280932

[27] SATO, K. and WATANABE, T. (2004). Moments of last exist times for Lévy processes. *Ann. Inst. H. Poincaré Probab. Statist.* **40** 207–225. MR2044816

[28] STROOCK, D. W. (1994). *A Concise Introduction to the Theory of Integration*, 2nd. ed. Birkhäuser, Boston. MR1267228

[29] STOECKL, A. and WAKOLBINGER, A. (1994). On clan-recurrence and -transience in time stationary branching Brownian particle systems. In *Measure-Valued Processes, Stochastic Partial Differential Equations, and Interacting Systems* (D. A. Dawson, ed.). *CRM Proc. Lecture Notes* **5** 213–219. Amer. Math. Soc., Providence, RI. MR1278296

[30] TALARCZYK, A. (2007). A functional ergodic theorem for the occupation time process of a branching system. *Statist. Probab. Lett.* To appear.
T. BOJDECKI  
A. TALARCZYK  
INSTITUTE OF MATHEMATICS  
UNIVERSITY OF WARSAW  
UL. BANACHA 2  
02-097 WARSAW  
POLAND  
E-MAIL: tobojd@mimuw.edu.pl  
  annatal@mimuw.edu.pl  

L. G. GOROSTIZA  
DEPARTMENT OF MATHEMATICS  
CENTRO DE INVESTIGACIÓN  
  Y DE ESTUDIOS AVANZADOS  
A.P. 14-740  
MÉXICO 07000 D.F.  
MEXICO  
E-MAIL: lgorosti@math.cinvestav.mx